\documentclass{amsart}

\usepackage{amssymb,latexsym}
\usepackage{palatino}
\usepackage{mathrsfs}

\newcommand{\Hom}{\operatorname{Hom}}

\newcommand{\tensor}{\otimes}

\newcommand{\congruent}{\equiv}
\newcommand{\iso}{\simeq}

\newcommand{\Int}{\operatorname{Int}}

\newcommand\R{\mathbf{R}}    
\newcommand\C{\mathbf{C}}    
\newcommand\Z{\mathbf{Z}}    
\newcommand\G{\mathbf{G}}    

\newcommand\TTT{\mathbf{T}}

\newcommand\AAA{\mathscr{A}}
\newcommand\BBB{\mathscr{B}}

\newcommand\aaa{\mathbf{a}}
\newcommand\bbb{\mathbf{b}}

\newcommand\e{\varepsilon}

\newcommand{\XX}{\mathscr{X}}
\newcommand{\lie}[1]{\mathfrak{#1}}
\newcommand{\glie}{\lie{g}}
\newcommand{\plie}{\lie{p}}

\newcommand{\UU}{\mathcal{U}}
\newcommand{\NN}{\mathcal{N}}

\newcommand{\Lie}{\operatorname{Lie}}
\newcommand{\Ad}{\operatorname{Ad}}

\newtheorem{theorem}{Theorem}
\newtheorem{prop}[theorem]{Proposition}
\newtheorem{cor}[theorem]{Corollary}
\newtheorem{lem}[theorem]{Lemma}
\theoremstyle{remark}
\newtheorem{rem}[theorem]{Remark}

\begin{document}
\bibliographystyle{hamsalpha} 
\author{George J. McNinch} 
\author{Eric Sommers} 

\dedicatory{To Robert Steinberg, on his 80th birthday.}

\thanks{Both authors were supported in part  
by the National Science Foundation;
the first author, by DMS-9970301 and the second author, by DMS-0070674.}

\title{Component groups of unipotent centralizers in good characteristic}
\address[George McNinch]{Department of Mathematics \\
         Room 255 Hurley Building \\
         University of Notre Dame \\
         Notre Dame, Indiana 46556-5683 \\
         USA}
\email{mcninch.1@nd.edu}

\address[Eric Sommers]{Department of Mathematics and Statistics\\
         University of Massachusetts at Amherst \\
         Amherst, MA 01003 \\
         USA}
\email{esommers@math.umass.edu}

\date{August 5, 2002}

\begin{abstract}
  Let $G$ be a connected, reductive group over an algebraically closed
  field of good characteristic.  For $u \in G$ unipotent, we describe
  the conjugacy classes in the component group $A(u)$ of the
  centralizer of $u$.  Our results extend work of the second author
  done for simple, adjoint $G$ over the complex numbers.
  
  When $G$ is simple and adjoint, the previous work of the second author makes
  our description combinatorial and explicit; moreover, it turns out
  that knowledge of the conjugacy classes suffices to determine the
  group structure of $A(u)$. Thus we obtain the result, previously
  known through case-checking, that the structure of the component
  group $A(u)$ is independent of good characteristic.
\end{abstract}

\maketitle

Throughout this note, $G$ will denote a connected and reductive
algebraic group $G$ over the algebraically closed field $k$. For the
most part, the characteristic $p \ge 0$ of $k$ is assumed to be
\emph{good} for $G$ (see \S \ref{section:simple} for the definition).

The main objective of our note is to extend the work of the second
author \cite{sommers-generalized-bala-carter} describing the component
groups of unipotent (or nilpotent) centralizers. We recall a few
definitions before stating the main result.

A \emph{pseudo-Levi} subgroup $L$ of $G$ is the connected centralizer
$C_G^o(s)$ of a semisimple element $s \in G$. The reductive group $L$
contains a maximal torus $T$ of $G$, and so $L$ is generated by $T$
together with the 1 dimensional unipotent subgroups corresponding to a
subsystem $R_L$ of the root system $R$ of $G$; in \S
\ref{sec:semisimple-centralizers} we make explicit which subsystems
$R_L$ arise in this way when $G$ is quasisimple.

Let $u \in G$ be a unipotent element, and let $A(u) = C_G(u)/C_G^o(u)$
be the group of components (``component group'') of the centralizer of
$u$.  We are concerned with the structure of the group $A(u)$ (more
precisely: with its conjugacy classes).

Consider the set of all triples 
\begin{equation}
  \label{eq:triples}
  (L,tZ^o,u)
\end{equation}
where $L$ is a pseudo-Levi subgroup with center $Z = Z(L)$, the coset
$tZ^o \in Z/Z^o$ has the property that $L = C_G^o(tZ^o)$, and $u \in
L$ is a distinguished unipotent element.  
\begin{theorem}
  \label{theorem:main-theorem}
  Let $G$ be connected and reductive in good characteristic. The map
  $$(L,sZ^o,u) \mapsto (u,sC^o_G(u))$$
  yields a bijection between:
  $G$-conjugacy classes of triples as in \eqref{eq:triples}, and
  $G$-conjugacy classes of pairs $(u,x)$ where $u \in G$ is unipotent
  and $x$ is an element in $A(u)$.
\end{theorem}

The theorem is proved, after some preliminaries, in \S\ref{sec:main-proof}.

\begin{rem}
  The $G$-conjugacy classes of pairs $(u,x)$ as in the statement of
  the theorem are in obvious bijection with $G$-conjugacy classes of
  pairs $(u,C)$ where $u \in G$ is unipotent and $C \subset A(u)$ is a
  conjugacy class.
\end{rem}

\begin{rem}
  Assume that $G$ is simple and adjoint.  We show in \S
  \ref{sec:explicit} that our work indeed extends the results of the
  second author.  If $u \in G$ is unipotent, we find as a consequence
  of Theorem \ref{theorem:main-theorem} that the conjugacy classes in
  $A(u)$ are in bijection with $C_G(u)$-conjugacy classes of
  pseudo-Levi subgroups $L$ containing $u$ as a distinguished
  unipotent element; this was proved for $k=\C$ in
  \cite{sommers-generalized-bala-carter}. It follows that $A(u) \iso
  A(\hat u)$ where $\hat u$ is a unipotent element in the
  corresponding group over $\C$ with the same labeled diagram as $u$.
  This isomorphism was known previously by case-checking in the
  exceptional groups; see especially \cite{mizuno}.  The structure of
  $A(u)$ for the exceptional groups when $k=\C$ is originally due to
  Alekseevski \cite{alekseevski}.
\end{rem}

\begin{rem}
  Our proof of Theorem \ref{theorem:main-theorem} is free of
  case-checking, with the following caveats. We use Pommerening's
  proof of the Bala-Carter theorem (specifically, we use the
  construction of ``associated co-characters'' for unipotent elements)
  in the proof of Proposition \ref{prop:connected-unipotent}.
  Moreover, we use work of Premet to find Levi factors in the
  centralizer of a unipotent element; see Proposition
  \ref{prop:levi-decomposition}.
\end{rem}
 
The authors would like to thank the referee for pointing out an
oversight and suggesting the use of Jantzen's result (Proposition
\ref{prop:weak-mostow}) to prove Proposition \ref{prop:assoc-in-PL}.
Upon completion of this paper, 
we learned that Premet has also given a case-free proof of
Theorem \ref{theorem:adjoint}.

\section{Reductive algebraic groups}
\label{section:simple}

Fix $T \subset B \subset G$, where $T$ is a maximal torus and $B$ a
Borel subgroup.  Let $(X,R,Y,R^\vee)$ denote the root datum of the
reductive group $G$ with respect to $T$; thus $X = X^*(T)$ is the
character group, and $R \subset X$ is the set of weights of $T$ on
$\glie$. Fix $S \subset R$ a system of simple roots.  

When $R$ is irreducible, the root with
maximal height (with respect to $S$) will be denoted $\tilde \alpha$.
Write
\begin{equation}
  \label{eq:high-root}
  \tilde \alpha = \sum_{\beta \in S} a_\beta \beta
\end{equation}
for positive integers $a_\beta$.  The characteristic $p$ of $k$ is
said to be good for $G$ (or for $R$) if $p$ does not divide any
$a_\beta$.  So $p=0$ is good, and we may
simply list the bad (i.e. not good) primes: $p=2$ is bad unless $R =
A_r$, $p=3$ is bad if $R = G_2,F_4,E_r$, and $p=5$ is bad if $R=E_8$.

The prime $p$ is good for a general $R$ just in case it is good for
each irreducible component of $R$.

For a root $\alpha \in R$, let $\G_a \iso \XX_\alpha \subset G$ be the
corresponding root subgroup.

\section{Springer's isomorphism}

Let $\UU \subset G$ and $\NN \subset \glie$ denote respectively the
unipotent and nilpotent varieties. In characteristic 0, the
exponential is a $G$-equivariant isomorphism $\NN \to \UU$; in
good characteristic, one has the following substitute for the
exponential:

\begin{prop}
  \label{prop:springer-iso}
  There is a $G$-equivariant homeomorphism $\e:\NN \to \UU$.
  Moreover, if $R$ has no component of type $A_r$ for which $r \equiv
  -1 \pmod p$, there is such an $\e$ which is an isomorphism of
  varieties.
\end{prop}

\begin{proof}
  There is an isogeny $\pi:\tilde G \to G$ where $\tilde G = \prod_i
  G_i \times T$ with $T$ a torus and each $G_i$ a simply connected,
  quasisimple group; see e.g. \cite[Theorem
  9.6.5]{springer-LAG}.  Let $\tilde \NN$ and
  $\tilde \UU$ denote the corresponding varieties for $\tilde G$.
  Since the characteristic is good, it has been proved by Springer
  \cite{springer-unipotent-iso} that there is a $\tilde G$-equivariant
  isomorphism of varieties $\tilde \e:\tilde \NN \to \tilde \UU$; for
  another proof, see \cite{Bard-Rich-LunaSlice}.
  
  It follows from \cite[Lemma 24]{Mc:sub-principal} that $\pi$
  restricts to a homeomorphism $\pi_{\mid \tilde \UU}:\tilde \UU \to
  \UU$, and that $d\pi$ restricts to a homeomorphism $d\pi_{\mid
    \tilde \NN}:\tilde \NN \to \NN$.  Since the characteristic is
  good, $d\pi$ is bijective provided that $R \not = A_r$ when $r
  \equiv -1 \pmod p$; see the summary in \cite[0.13]{hum-conjugacy}.
  It follows from the remaining assertion in \cite[Lemma
  24]{Mc:sub-principal} that $\pi_{\mid \tilde \UU}$ and $d\pi_{\mid
    \tilde \NN}$ are isomorphisms of varieties when $d\pi$ is
  bijective, whence the proposition.
\end{proof}

In what follows, \emph{we fix an equivariant homeomorphism 
$\e:\NN \to \UU$, to which we will refer 
without further comment.}

\section{Associated co-characters}

Recall that a unipotent $u \in G$ is said to be \emph{distinguished} 
if the connected center $Z^o(G)$ of $G$ is a maximal torus of 
$C_G(u)$. A nilpotent element $X \in \glie$ is then distinguished if $\e(X)$
has that property.

Let $X \in \glie$ be nilpotent.  If $X$ is not distinguished, there is
a Levi subgroup $L$ of $G$ for which $X \in \Lie(L)$ is distinguished.

A co-character
$\phi:k^\times \to G$ is said to be associated to $X$ if
\begin{equation*}
  \Ad \phi(t) X = t^2 X \quad \text{for each} \quad t \in k^\times,
\end{equation*}
and if the image of $\phi$ lies in the 
derived group of some Levi subgroup $L$ for which $X \in \Lie(L)$ is
distinguished.

A co-character $\phi$ is  associated to a unipotent $u \in G$
if it is associated to $X=\e^{-1}(u)$.

\begin{prop}
  \label{prop:co-characters}
  Let $u \in G$ be unipotent. Then there exist co-characters
  associated to $u$, and any two such are conjugate by an element of
  $C_G^o(u)$.
\end{prop}

\begin{proof}
  This is proved in \cite[Lemma 5.3]{jantzen:Nilpotent}.
\end{proof}

\begin{rem}
  The existence of associated co-characters asserted in
  the the previous proposition relies in an essential way on
  Pommerening's proof \cite{Pommerening} of the Bala-Carter theorem in
  good characteristic.
\end{rem}

Let $\phi$ be a co-character associated to the unipotent $u\in G$, and
let $\glie(i)$ be the $i$-weight space for $\Ad \circ \phi(k^\times)$, $i \in
\Z$.  Let $\plie = \bigoplus_{i \ge 0} \glie(i)$. Then $\plie =
\Lie(P)$ for a parabolic subgroup $P$ of $G$; $P$ is known as the canonical
parabolic associated with $u$.
\begin{prop}
  \label{prop:canonical-parabolic}
  Let $u \in G$ be unipotent.
  The parabolic subgroup $P$ is independent of the choice of
  associated co-character $\phi$ for $u$. Moreover,
  $C_G(u) \le P$.
\end{prop}
\begin{proof}
   \cite[Prop. 5.9]{jantzen:Nilpotent}
\end{proof}

\begin{rem}
  The proof that $C_G(u) \subset P$ is somewhat subtle in positive
  characteristic. Let $X = \e^{-1}(u)$. In characteristic 0, the
  assertion $C_G^o(u) \subset P$ is a consequence of the Lie algebra
  analogue $\lie{c}_\glie(X) \subset \lie{p}$ which follows from the
  Jacobson-Morozov theorem.  (The fact that the full centralizer lies
  in $P$ is then a consequence of the unicity of the canonical
  parabolic $P$). In good characteristic, the required assertion for
  the Lie algebra was proved by Spaltenstein, and independently by
  Premet; see the references in \cite[\S 5]{jantzen:Nilpotent}.  In
  the positive characteristic case, the transition to the group is
  more subtle; again see \emph{loc. cit.}
\end{rem}

\section{The Levi decomposition of a unipotent centralizer}

In characteristic $p>0$, a linear algebraic group can fail to have a
Levi decomposition. Moreover, even when they exist, two Levi factors
need not in general be conjugate. If $u \in G$ is unipotent and
the characteristic is good for $G$, the connected centralizer
$C_G^o(u)$ does have  a Levi decomposition, thanks to work of Premet.
More precisely:

\begin{prop}
  \label{prop:levi-decomposition}
  Let $u \in G$ be unipotent, let $P$ be the canonical parabolic
  associated with $u$ (see Proposition \ref{prop:canonical-parabolic}),
  and let $U_P$ be the unipotent radical of $P$.
  \begin{enumerate}
  \item $R(u) = C_G(u) \cap U_P$ is the unipotent radical of $C_G(u)$.
  \item For any co-character $\phi$ associated with $u$, the
    centralizer $C_\phi$ of $\phi$ in $C_G(u)$ is a Levi factor of
    $C_G(u)$; i.e. $C_\phi$ is reductive and $C_G(u) = C_\phi \cdot
    R(u)$.
  \item If $\phi,\phi'$ are two co-characters associated to $u$, then
    $C_\phi$ and $C_{\phi'}$ are conjugate by an element in
    $C_G^o(u)$.
  \end{enumerate}
\end{prop}

\begin{proof}
  \cite[\S 5.10, 5.11]{jantzen:Nilpotent}.
\end{proof}

\begin{rem}
  The proof that $R(u)$ is a connected (normal, unipotent) group is
  elementary, as is the fact that $C_G(u) = C_\phi \cdot R(u)$.  The
  proof that $C_\phi$ is reductive depends on work of Premet, and
  ultimately involves case-checking in small characteristics for
  exceptional groups.
\end{rem}

\section{Semisimple representatives}

If $H$ is a linear algebraic group, in characteristic
0 one may always represent a coset $tH^o \in H/H^o$ by 
a semisimple element $t \in H$. In characteristic $p>0$ this
is no longer true in general (e.g. if $[H:H^o] \congruent 0 \pmod{p}$).

Let now $G$ be connected, reductive in good characteristic and suppose
$u \in G$ is unipotent.  Despite the above difficulty, we may always
choose semisimple representatives for the elements in the component
group $A(u)$.

\begin{prop}
  \label{prop:connected-unipotent}
  Let $u \in G$ be unipotent, and suppose $v \in C_G(u)$ is
  also unipotent. Then $v \in C_G^o(u)$.
\end{prop}

\begin{proof}
  The proposition follows from \cite[III.3.15]{springer-steinberg}.
  Note that in \emph{loc. cit.} $G$ is assumed semisimple, but the
  argument works for all reductive $G$ in view of Proposition
  \ref{prop:springer-iso}.
\end{proof}

\begin{cor}
  \label{cor:semisimple-representative}
  Let $u \in G$ be unipotent.  Then each element of the component
  group $A(u)$ may be represented by a semisimple element $s \in
  C_G(u)$.
\end{cor}

\begin{proof}
  Let $g \in C_G(u)$, and let $g=g_s g_u$ be its Jordan decomposition
  where $g_s$ is semisimple and $g_u$ is unipotent.  Proposition
  \ref{prop:connected-unipotent} implies that $g_u \in C_G^o(u)$,
  whence the corollary.
\end{proof}

\section{Pseudo-Levi subgroups}

We collect here a few results on pseudo-Levi subgroups which will be
needed in the proof of Theorem \ref{theorem:main-theorem}.  Recall
that by a pseudo-Levi subgroup, we mean the connected centralizer of a
semisimple element of $G$.

\begin{lem}
  \label{lem:ss-centralizer}
  Let $S \subset T$ be a subset. Then
  $C_G^o(S)$ is a reductive subgroup of $G$, and is generated by $T$
  together with the root subgroups $\XX_\alpha$ for which $\alpha(S) =
  1$.
\end{lem}

\begin{proof}
  \cite[II \S 4.1]{springer-steinberg}.
\end{proof}

\begin{prop}
  \label{prop:misc-plevi}
  Let $L=C_G^o(t)$ with $t \in G$ semisimple. Write $Z$ for the
  center of $L$.
  \begin{enumerate}
  \item $L = C_G^o(tZ^o)$.
  \item Let $S$ be a torus in $C_G^o(t)$, and let $M=C_G^o(tS)$.
    There is a non-empty open subset $U \subseteq tS$ such that
    $M = C_G^o(x)$ for each  $x \in U$. In particular, $M$ is again
    a pseudo-Levi subgroup of $G$.  If $Z_1$ denotes the center of
    $M$, then $M = C_G^o(tZ_1^o)$.

  \item There is a non-empty open subset $U
    \subseteq tZ^o$ such that $L = C_G^o(x)$ for each $x \in U$.
  \end{enumerate}
\end{prop}

\begin{proof}
  (1) is straightforward to verify. 
  
For (2), we may suppose that $t$ and $S$ are in $T$. Let
  $R' =\{\alpha \in R \mid \alpha(tS) = 1\}$. Then $R' \subseteq R_x =
  \{\alpha \in R \mid \alpha(x) = 1\}$ for any $x \in tS$.  Since $tS$
  is an irreducible variety, the intersection of non-empty open
  subsets
  \begin{equation*}
    U=\bigcap_{\alpha \in R \setminus R'} \{x \in tS \mid \alpha(x) \not = 1\}
  \end{equation*}
  is itself open and non-empty; moreover, it is clear that $R_x = R'$
  whenever $x \in U$, so the first assertion of (2) follows from Lemma
  \ref{lem:ss-centralizer}.
  
  For the final assertion of (2), first note that $M = C_G^o(t,S) =
  C_L^o(S)$ is a Levi subgroup of $L$.  By \cite[Prop.
  1.21]{DigneMichel} we have $M = C_L^o(Z_1^o)$; since $t$ is central
  in $M$, we have also $M = C_L^o(tZ_1^o)$.  Since certainly
  $C_G^o(tZ_1^o) \subseteq C_G^o(t) = L$, we deduce that $M =
  C_G^o(tZ_1^o)$ as desired.

  (3) follows from (1) and (2) with $S = Z^o$.
\end{proof}

\begin{prop}
  \label{prop:pl-good-char}
  Let $G$ be connected and reductive.  If the characteristic $p$ of
  $k$ is good for $G$, and if $L$ is a pseudo-Levi subgroup of $G$,
  then $p$ is good for $L$ as well.
\end{prop}

\begin{proof}
  As in the proof of Proposition \ref{prop:springer-iso}, let
  $\pi:\tilde G \to G$ be an isogeny where $\tilde G = \prod_i G_i
  \times S$ with $S$ a torus and each $G_i$ a simply connected
  quasisimple group. Let $L = C_G^o(t)$. If $\pi(\tilde t) = t$ and
  $\tilde L = C^o_{\tilde G}(\tilde t)$, then Lemma
  \ref{lem:ss-centralizer} shows that $\pi(\tilde L) = L$. Since $p$
  is good for $L$ if and only if it is good for $\tilde L$, we may
  replace $G$ by $\tilde G$. Since $L = \prod_i (L \cap G_i) \times
  S$, it suffices to suppose that $G$ is quasisimple.
  
  According to \cite[\S4.1,4.3]{springer-steinberg} $p$ is good for
  $G$ if and only if $\Z R/ \Z R_1$ has no $p$-torsion for any
  (integrally) closed subsystem $R_1$ of $R$.  Since the root system
  $R_L$ of $L$ is one such subsystem, it readily follows that $p$ is good
  for any irreducible component of  $R_L$.
\end{proof}

\begin{lem}
  \label{lem:finite-rep}
  Let $L$ be a pseudo-Levi subgroup of $G$. Then $L=C_G^o(s)$ for a
  semisimple element $s \in G$ of finite order.
\end{lem}

\begin{proof}
  Let $Z$ denote the (full) center of $L$.  By \cite[Exerc.  3.2.10
  5(b)]{springer-LAG}, the elements of $Z$ which
  have finite order are dense in the diagonalizable group $Z$. Now
  choose $t \in Z$ such that $L=C_G^o(t)$, and let $U \subset tZ^o$ be
  an open set as in Proposition \ref{prop:misc-plevi}(3). Then $U$ is
  also open in $Z$ and hence contains an element $s$ of finite order.
\end{proof}

\section{Semisimple automorphisms of reductive groups}

If $H$ is any linear algebraic group, an automorphism $\sigma$ of $H$
is semisimple if there is a linear algebraic group $H'$ 
with $H \lhd H'$ such that $\sigma = \Int(x)_{\mid H}$ for some semisimple
$x \in H'$ (where $\Int(x)$ denotes the inner automorphism determined by $x$).

\begin{prop}
  \label{prop:Steinberg-ss}
  Let $H$ be a connected linear algebraic group, and let
  $\sigma$ be a semisimple automorphism of $H$.  
  Then $\sigma$ fixes a
    Borel subgroup $B$ of $H$, and a maximal torus $T \subset B$.
\end{prop}
\begin{proof}
  \cite[Theorem 7.5]{steinberg-endomorphisms}. 
\end{proof}

\begin{lem}
  \label{lem:product-decomposition}
  Let $A$ be a connected commutative linear algebraic group, let
  $\sigma$ be a semisimple automorphism of $A$, and let $A^\sigma$ be
  the fixed points of $\sigma$ on $A$.  Then each element $a \in A$
  can be written $$a=x \cdot \sigma(y)y^{-1}$$
  for $x \in A^\sigma$ and
  $y \in A$.
\end{lem}

\begin{proof}
  The homomorphism $\phi:A^\sigma \times A \to A$ given by $\phi(x,y)
  = x \cdot \sigma(y)y^{-1}$ has surjective differential by
  \cite[Corollary 5.4.5(ii)]{springer-LAG}, so
  $\phi$ is surjective.
\end{proof}

\begin{prop}
  \label{prop:key-conjugacy-result}
  Let $H$ be a reductive algebraic group, and suppose the images of
  the semisimple elements $t,t' \in H$ lie in the same conjugacy class
  in $H/H^o$. Then there is $g \in H$ and a semisimple $s \in C_H^o(t)$ such that
  $gt'g^{-1} = ts.$
\end{prop}

\begin{proof}
  Replacing $t'$ by $ht'h^{-1}$ for suitable $h \in H$, we can suppose
  that $t$ and $t'$ have the same image in $H/H^o$.
  
  Applying Proposition \ref{prop:Steinberg-ss} we can find $T \subset
  B$ where $T$ and $B$ are respectively an $\Int(t)$-stable maximal
  torus and Borel group.  Similarly, we can find an
  $\Int(t')$ stable  $T' \subset B'$.
  
  Choose $g \in H$ with $B = g^{-1}B'g$.  Then $g^{-1}T'g$ is a
  sub-torus of $B$. Replacing $g$ by $bg$ for some $b \in B$, we can
  arrange that $g^{-1}T'g = T$; replacing $t'$ by $gt'g^{-1}$, we see
  that  $T \subset B$ is both $\Int(t)$-stable and
  $\Int(t')$-stable.
  
  Thus $n=t^{-1}t'$ is in the normalizer of $T$ in $H^o$.  Since
  $\Int(n)$ fixes $B$, and since the Weyl group $N_{H^o}(T)/T$ acts
  simply transitively on the set of Borel subgroups containing $T$, we
  deduce that $n \in T$.  We can therefore write $t' = ta$ for $a \in
  T$.  By Lemma \ref{lem:product-decomposition} we can write $a = x
  t^{-1}yty^{-1}$ for some $x \in C_T(t)$ and $y \in T$.
  Let $g = ty^{-1}$. Then one readily checks that
  \begin{equation*}
    gt'g^{-1}  = tx
  \end{equation*}
  and the proof is complete.
\end{proof}

\begin{cor}
  \label{cor:conjugacy-result}
  Let $H$ be a linear algebraic group. Suppose that $\mathcal{M}$ is a
  collection of Levi factors of $H$ which are all conjugate under $H$.
  If the semisimple elements $t,t' \in H$ lie in the same conjugacy
  class in $H/H^o$, and if $t,t' \in \cup_{M \in \mathcal{M}} M$, then
  there is $g \in H$ and a semisimple element $s \in C_H^o(t)$ such
  that $gt'g^{-1} = ts$.
\end{cor}

\begin{proof}
  Choose $M,M' \in \mathcal{M}$ with $t \in M$ and $t' \in M'$.  Since
  $M$ and $M'$ are $H$-conjugate, replacing $t'$ by an $H$-conjugate
  permits us to suppose that $t,t' \in M$.  Since
  $M$ is reductive, we deduce the result from Proposition
  \ref{prop:key-conjugacy-result}.
\end{proof}

We require one further property of pseudo-Levi subgroups, whose proof
depends on Proposition \ref{prop:Steinberg-ss} and on the following
version of a result of Mostow recently obtained by
Jantzen \cite[11.24]{jantzen:Nilpotent}:
\begin{prop}
  \label{prop:weak-mostow}
  Let $\Gamma$ be an algebraic group which is a semidirect product of
  a (not necessarily connected) reductive group $M$ and a normal
  unipotent group $R$. Let $H \le \Gamma$ be a linearly reductive
  closed subgroup of $\Gamma$. Then there exists $r \in R$ with
  $rHr^{-1} \subset M$.
\end{prop}

\begin{prop}
  \label{prop:assoc-in-PL}
  Let $L$ be a pseudo-Levi subgroup of $G$ and $u \in L$ a
  distinguished unipotent element.  If a cocharacter of L is
  associated to $u$ in $L$, then that cocharacter is associated to $u$
  in $G$ as well.
\end{prop}

\begin{proof}
  Since all cocharacters associated to $u$ in $L$ are conjugate by
  $C_L^o(u)$ (Proposition \ref{prop:co-characters}), it suffices to
  find some cocharacter of $L$ which is associated to $u$ in both $L$
  and $G$.

  According to Lemma \ref{lem:finite-rep}, $L = C_G^o(s)$ for some
  semisimple element $s$ of finite order. The order of $s$ is then
  invertible in $k$, so the cyclic subgroup $H$ generated by $s$ is
  linearly reductive (all of its linear $k$-representations are
  completely reducible).
  
  Let $\phi$ be any cocharacter of $G$ associated to $u$, and consider
  the subgroup $N = \phi(k^\times) C_G(u)$ (i.e. the group generated by
  the centralizer, and by the image of $\phi$; this is the group
  $N(\e^{-1}(u))$ defined in \cite[2.10(2)]{jantzen:Nilpotent}).
  
  According to Proposition \ref{prop:levi-decomposition}, the
  centralizer $C_\phi$ of $\phi$ in $C_G(u)$ is a Levi factor of $C_G(u)$.
  Now $C_\phi'=\phi(k^\times)\cdot C_\phi$ is a Levi factor
  of $N$.  Moreover, the image of $\phi$ is central in $C_\phi'$.
  
  Now take $\Gamma = N$ in Proposition \ref{prop:weak-mostow}. Then $H
  = \langle s \rangle$ is a linearly reductive subgroup of $\Gamma$.
  So there is an element $r$ in the unipotent radical of $C_G(u)$ such
  that $rsr^{-1}$ is in $C_\phi'$.  But then $rsr^{-1}$ centralizes the
  image of $\phi$, so that $s$ centralizes the image of
  $\phi'=\Int(r^{-1}) \circ \phi$. Thus, $\phi'$ is a cocharacter of $L$.
  
  We claim that $\phi'$ is associated to $u$ in $L$.  Since the map
  $\e:\NN\to \UU$ is $G$-equivariant and thus restricts to a
  homeomorphism $\NN(L) \to \UU(L)$, we must see that $\phi'$ is
  associated to $\e^{-1}(u)$. Thus, we only must verify that $\phi'$
  takes values in the derived group of $L$.
  
  For each maximal torus $S$ of $C_{\phi'}$, $u$ is distinguished in
  $M = C_G(S)$ and the image of $\phi'$ lies in the derived group
  $(M,M)$ (to see this last assertion, note that it holds for
  \emph{some} such $S$ since $\phi'$ is associated to $u$ in $G$, and
  hence for all such $S$ by conjugacy of maximal tori in $C_{\phi'}$).
  
  We may choose a maximal torus $S \le C_{\phi'}$ containing the
  connected center of $L$. Since $C_{\phi'}$ is normalized by $s$, we
  may also suppose by Proposition \ref{prop:Steinberg-ss} that $S$ is
  normalized by $s$. Then $C_S^o(s)$ is a torus in $C_L(u)$; since $u$
  is distinguished in $L$, we see that $C_S^o(s)$ is contained in (and
  hence coincides with) the connected center of $L$.
  
  The subgroup $M$ is normalized by $s$, and the proposition will
  follow if we can show that $\phi'$ takes values in the derived
  group of $C_M^o(s)$ (since $C_M^o(s) \subset L$).  
  
  We first claim that $C_S^o(s)$ is the maximal central torus of
  $C^o_M(s)$.  Indeed, if $C_S^o(s) \subset S'$ with $S'$ a central
  torus of $C^o_M(s)$, then $S'$ centralizes $s$ and $u$ so that $S'
  \subset C_L(u)$; since $C_S^o(s)$ is the unique maximal torus of
  $C_L(u)$, $S' = C_S^o(s)$ as claimed.

  The proposition is now a consequence of the lemma which follows.
\end{proof}

\begin{lem}
  Let $M$ be a connected, reductive group, and suppose that $\sigma$
  is a semisimple automorphism of $M$. If $S$ is a $\sigma$-stable
  central torus in $M$ and $C_S^o(\sigma)$ is the maximal central
  torus of $C^o_M(\sigma)$, then
\begin{equation}
  \label{eq:derived-=}
  (C_M^o(\sigma),C_M^o(\sigma)) =  C_{(M,M)}^o(\sigma).
\end{equation}
\end{lem}

\begin{proof}
  The inclusion
  $$(C_M^o(\sigma),C_M^o(\sigma)) \subseteq C_{(M,M)}^o(\sigma)$$
  is
  immediate (by \cite[2.2.8]{springer-LAG} the group
  on the left is connected; it is also evidently a $\sigma$-stable
  subgroup generated by commutators in $M$).
  
  On the other hand, according to \cite[9.4]{steinberg-endomorphisms},
  $N=C_{(M,M)}^o(\sigma)$ is reductive.  We claim that $N$ is
  semisimple; if that is so then $N=(N,N)$; since $N \subseteq
  C_M^o(\sigma)$, equality in \eqref{eq:derived-=} will follow.
  
  Write $Z$ for the connected center of $M$. Then $Z \cap (M,M)$ is
  finite; see e.g. \cite[8.1.6]{springer-LAG}.  Since $S \subseteq Z$,
  we see that $C_S^o(\sigma) \cap N$ is finite as well.
  
  Now let $T$ be any $\sigma$-stable maximal torus of $M$.  We know
  that $\Lie(M)$ is the sum of $\Lie((M,M))$ and $\Lie(T)$, since
  $\Lie((M,M))$ contains each non-zero $T$-weight space of $\Lie(M)$.
  It follows from \cite[Lemma 4.4.12]{springer-LAG} that the
  differential at $(1,1)$ of the product map $\mu:T \times (M,M) \to
  M$ is surjective.  Since $d\sigma$ is diagonalizable, $d\mu_{(1,1)}$
  restricts to a surjective map on $d\sigma$-eigenspaces (for each
  eigenvalue); especially, it restricts to a surjective map on the
  fixed points of $d\sigma$.  Reinterpreting the $d\sigma$-fixed
  points via \cite[5.4.4]{springer-LAG}, we see that the restriction
  of $d\mu_{(1,1)}$ to $\Lie(C_T(\sigma)) \oplus
  \Lie(C_{(M,M)}(\sigma))$ surjects onto $\Lie(C_M(\sigma))$.  It
  follows that $\mu$ restricts to a dominant morphism $\tilde
  \mu:C^o_T(\sigma) \times N \to C^o_M(\sigma)$; cf.
  \cite[4.3.6]{springer-LAG}.  Since $C^o_T(\sigma)$ normalizes $N$,
  the image is a subgroup.  As $C^o_M(\sigma)$ is connected, $\tilde
  \mu$ is surjective.
  
  Let $R$ denote the radical of $N$ ($R$ is the maximal central torus
  of $N$). By Proposition \ref{prop:Steinberg-ss}, $R$ is
  contained in $C_T(\sigma)$ for some $\sigma$-stable
  maximal torus $T$ of $M$. 
  Applying the considerations of the previous paragraph to this $T$,
  we get that $C^o_M(\sigma) = C^o_T(\sigma) \cdot N$. It follows that
  $R$ is a central torus in $C^o_M(\sigma)$.  Since $C_S^o(\sigma)$ is
  by assumption the maximal central torus of $C^o_M(\sigma)$, we have
  that $R \subseteq C_S^o(\sigma) \cap N$ is finite, so $R=1$ and $N$
  is indeed semisimple.
\end{proof}

\begin{rem}
  Though we shall not have occasion to use it here, the conclusion of
  Proposition \ref{prop:assoc-in-PL} is true more generally: $(*)$ if
  $L$ is a pseudo-Levi subgroup, and if $\phi$ is a cocharacter of $L$
  associated to a unipotent $u \in L$, then $\phi$ is associated to
  $u$ in $G$. This follows from Proposition \ref{prop:assoc-in-PL}
  together with the observation that a Levi subgroup of $L$ is a
  pseudo-Levi subgroup of $G$ (Proposition \ref{prop:misc-plevi}(2)),
  and that $(*)$ holds when $L$ is a Levi subgroup.
\end{rem}

\section{Establishing the main result}
\label{sec:main-proof}

Let $\AAA$ be the set of triples $\aaa=(L,tZ^o,u)$ where $L$ is a
pseudo-Levi subgroup of $G$ with center $Z$, $tZ^o \in Z/Z^o$
satisfies $C_G^o(tZ^o) = L$, and $u \in L$ is a distinguished
unipotent element. The action of $G$ on itself by inner automorphisms
determines an action of $G$ on $\AAA$.

For $\aaa = (L,tZ^o,u) \in \AAA$, we set $u(\aaa) = u$, and we write
$c(\aaa) \subset A(u)$ for the element $c(\aaa) = tC_G^o(u)$.

Let $\BBB$ be the set of all pairs $(u,x)$ where $u \in G$ is
unipotent and $x \in A(u)$. The action of $G$ on itself by inner
automorphisms yields an action of $G$ on $\BBB$.

To $\aaa \in \AAA$ we associate the pair $\Phi(\aaa) =
(u(\aaa),c(\aaa)) \in \BBB$.

\begin{lem}
  \label{lem:surjective}
  Let $(u,c) \in \BBB$. Then there is $\aaa \in \AAA$ with 
  $\Phi(\aaa) = (u,c)$.
\end{lem}

\begin{proof}
  Choose a semisimple $t \in C_G(u)$ whose image in $A(u)$ represents
  $c$ (Corollary
  \ref{cor:semisimple-representative}). Let $S$ be a maximal torus of
  $C_G^o(u,t)$.  Then $L = C_G^o(t,S) = C_G^o(tS)$ is a pseudo-Levi
  subgroup of $G$ containing $u$, and $L = C_G^o(tZ^o)$ where $Z$
  denotes the center of $L$ (Proposition \ref{prop:misc-plevi}(2)).
  
  It remains to show that $\aaa = (L,tZ^o,u)$ is in $\AAA$, i.e. that
  $u$ is distinguished in $L$.  Let $A$ be a maximal torus of
  $C_L(u)$; we must show that $A$ is in the center of $L$.  Note that
  $A$ is a subtorus of $B=C_G^o(u,t)$ and that $A$ centralizes
  $S$. In particular, $A$ is contained in the Cartan subgroup
  $H=C_B(S)$; by \cite[Prop 6.4.2]{springer-LAG} $H$
  is nilpotent and $S$ is its unique maximal torus.  Thus $A$ is
  contained in $S$, hence $A$ is central in $L$.
  
  It is clear that $\Phi(\aaa) = (u,c)$; this completes the proof.
\end{proof}

\begin{lem}
  \label{lem:injective}
  Let $\aaa,\bbb \in \AAA$, and suppose that $u = u(\aaa) = u(\bbb)$.
  If $c(\aaa)$ and $c(\bbb)$ are conjugate in $A(u)$, then 
  there is $g \in C_G(u)$ with $\aaa = g\bbb$.
\end{lem}

\begin{proof}
  Write $\aaa = (L,tZ^o,u)$ and $\bbb = (L',t'{Z'}^o,u)$.  By
  Proposition \ref{prop:misc-plevi}(3), we may choose the
  representatives $t,t'$ such that $L = C_G^o(t)$ and $L'=C_G^o(t')$.
  
  Let $\phi:k^\times \to L$ be a co-character associated to $u$ for
  the pseudo-Levi subgroup $L$; see Propositions
  \ref{prop:pl-good-char} and \ref{prop:co-characters}.  Then $\phi$
  is associated to $u$ in $G$ as well; see Proposition
  \ref{prop:assoc-in-PL}.  Evidently, $t \in C_\phi$, where $C_\phi$
  is the Levi factor of $C_G(u)$ of Proposition
  \ref{prop:levi-decomposition}. Similarly, $t'$ lies in a Levi factor
  $C_{\phi'}$ of $C_G(u)$.
  
  Consider the collection $\mathcal{M} = \{C_\phi \mid \phi$ is
  associated to $u\}$ of Levi factors of $C_G(u)$.  Then any two Levi
  factors in $\mathcal{M}$ are conjugate under $C_G^o(u)$ by
  Proposition \ref{prop:levi-decomposition}. The previous paragraph
  shows that $t,t' \in \cup_{M \in \mathcal{M}} M$, hence we may apply
  Corollary \ref{cor:conjugacy-result}. That corollary yields $g\in
  C_G(u)$ and a semisimple $s \in C_G^o(u,t)$, such that $gt'g^{-1} =
  ts$.
  
  Choose a maximal torus $S$ of $C_G^o(u,t)$ containing $s$.
  Then $S \subset L$ and $S$ centralizes $u$; since $u$ is distinguished in
  $L$, it follows that $s \in S \subset Z^o$.
  We have
  \begin{equation*}
    gL'g^{-1} = C_G^o(gt'g^{-1}) = C_G^o(ts).
  \end{equation*}
  Since $s \in Z^o$, we find that $L \subseteq C_G^o(ts)$. Thus $\dim
  L' \ge \dim L$. A symmetric argument shows that $\dim L' \le \dim
  L$, hence equality holds. We deduce that $gL'g^{-1} = C_G^o(ts) = L$,
  and so $g\bbb = \aaa$ as desired.
\end{proof}

\begin{proof}[Proof of Theorem \ref{theorem:main-theorem}]
  In the notation introduced in this section, Theorem
  \ref{theorem:main-theorem} is equivalent to: $\Phi$ induces a
  bijection from the set $\AAA/G$ of $G$-orbits on $\AAA$ to the set
  $\BBB/G$ of $G$-orbits on $\BBB$.
  
  First note that $\Phi(g\aaa) = g\Phi(\aaa)$ for each $\aaa \in
  \AAA$, so that indeed $\Phi$ induces a well-defined map 
  $\overline{\Phi}:\AAA/G \to \BBB/G$.
  Lemma \ref{lem:surjective} implies that $\Phi$ itself is
  surjective, hence also $\overline{\Phi}$ is surjective. Finally,
  Lemma \ref{lem:injective} shows that that
  $\overline{\Phi}$ is injective; this proves the theorem.
\end{proof}

\section{Centralizers of semisimple elements in quasisimple groups}

\label{sec:semisimple-centralizers}

In this section, we characterize the pseudo-Levi subgroups of $G$ when
the root system is irreducible (i.e. when $G$ is quasisimple); the
results are applied in the next section. The characterization we give
is well-known (certainly in characteristic 0), but as we have not
located an adequate reference (see Remark \ref{rem:deriziotis} below),
and since the arguments are not too lengthy, we include most details.

Let $T$ be any torus over $k$ with co-character group $Y$ (in the
application, we take $T$ to be a maximal torus of $G$). We denote by
$V = Y \tensor \R$ the extension of $Y$ to a real vector space, and by
$\TTT = V/Y$ the resulting compact (topological) torus.  If $X$ is the
character group of $T$, then $X$ identifies naturally with the
Pontrjagin dual $\hat \TTT = \Hom(\TTT,\R/\Z)$ of $\TTT$ [note that we
regard $\R/\Z$ as a \emph{multiplicative} group].  The
following lemma due to T. A.  Springer may be found in
\cite[\S5.1]{steinberg-endomorphisms}
\begin{lem}
  \label{lem:springer-dual-lemma}
  \begin{enumerate}
  \item[(a)] For each $t \in T$, there is $t' \in \TTT$ with the
    following property: 
    \begin{equation*}
      (*) \quad \text{for each $\lambda \in X$, $\lambda(t) = 1$ if
    and only if $\lambda(t') = 1$.}
    \end{equation*}   
  \item[(b)]  Conversely, if $t' \in \TTT$ has finite order,
    relatively prime to $p$ if $p>0$, there is $t \in T$
    for which $(*)$ holds.
  \end{enumerate}
\end{lem}

Unless explicitly stated otherwise,\emph{ we suppose in this section that
$R$ is irreducible, so that $G$ is quasisimple.}

Let $\tilde S = S \cup \{\alpha_0\}$ where $\alpha_0 = -\tilde
\alpha$; thus $\tilde S$ labels the vertices of the extended Dynkin
diagram of the root system $R$. For any subset $J \subsetneq \tilde
S$, let $R_J = \Z J \cap R$.  
Note that we do not require the characteristic to be good for $G$ in
this section.
\begin{lem}
  \label{lem:R_t=R_J}
  Let $T$ be our fixed maximal torus of $G$, and let $\TTT$ be the
  corresponding compact topological torus.  For $t \in \TTT$, put $R_t
  = \{\alpha \in R \mid \alpha(t) = 1\}$.  Then there is $J \subsetneq
  \tilde S$ such that $R_t$ is conjugate to $R_J$ by an element of
  $W$, the Weyl group of $R$.
\end{lem}

\begin{proof}
  Let $\tilde t \in V$ represent $t \in \TTT$. For some element
  $\tilde w$ of the affine Weyl group $W_a = W \cdot \Z Y$, $\tilde
  w\tilde t$ lies in the fundamental alcove $A_o$ in $V$ (whose walls
  are labeled by $\tilde S$).  The image of $\tilde w \tilde t$ in $V$
  is then $wt$, where $w$ is the image of $\tilde w$ in the finite
  Weyl group $W$, and $R_{wt} = w^{-1}R_t$. Thus, we suppose that $t$
  can be represented by a vector in $A_o$. In that case, let $J =
  \{\alpha \in \tilde S \mid \alpha(t) = 1\}$. Then the equality of $R_t$
  and $R_J$ is proved in \cite[Lemma 5.4]{lusztig-unipotent-p-adic}
  (in \emph{loc.  cit.}, Lusztig works instead with the complex torus
  $Y \tensor \C / Y$, but his argument is readily adapted to the
  current situation).
\end{proof}

For a subset $J \subsetneq \tilde S$, we consider the subgroup
\begin{equation*}
  L_J = \langle T, \XX_\alpha \mid \alpha \in R_J \rangle.
\end{equation*}

\begin{prop}
  \label{prop:semisimple-centralizer}
  Let $t \in G$ be semisimple. Then $C_G^o(t)$ is conjugate to a
  subgroup $L_J$ for some $J \subsetneq \tilde S$.  
\end{prop}

\begin{proof}
  We may suppose that $t \in T$. Set $R_t = \{\alpha \in R \mid
  \alpha(t) = 1\}$. According to Lemma \ref{lem:ss-centralizer},
  $C_G^o(t)$ is generated by $T$ and the $\XX_\alpha$ with $\alpha \in
  R_t$. With notations as before, choose $t' \in \TTT$ with the
  property $(*)$ of Lemma \ref{lem:springer-dual-lemma} for $t$. Thus
  $R_t = R_{t'}$. 
  Lemma \ref{lem:R_t=R_J} implies that $R_t$ and $R_J$ are $W$-conjugate
  for some  $J \subsetneq \tilde S$; thus
  $C_G^o(t)$ is conjugate in $G$ to $L_J$ as desired.
\end{proof}

\begin{rem}
  \label{rem:deriziotis}
  When $k$ is an algebraic closure of a finite field, Proposition
  \ref{prop:semisimple-centralizer} was proved by D. I. Deriziotis,
  and is stated in \cite[2.15]{hum-conjugacy}. See the last paragraph
  of \emph{loc. cit.} \S 2.15 for a discussion.
\end{rem}

In good characteristic, the converse of the previous proposition is
true as well:
\begin{prop}
  \label{prop:simple-pseudo-levi}
  Suppose that the characteristic of $k$ is good for $G$.
  Let $J \subsetneq \tilde S$, and let $Z$ be the center of $L_J$.
  There is $t \in Z$ with $L_J = C_G^o(t)$.
\end{prop}

\begin{proof}
  It suffices to suppose that $G$ is adjoint. In that case, there are
  vectors $\varpi_\alpha^\vee \in Y$, $\alpha \in S$, dual to the
  basis $S$ of $X$.  We suppose that $\alpha_0 \in J$, since otherwise
  $L_J$ is a Levi subgroup and the result holds (in all
  characteristics) e.g.  by
  \cite[6.4.3]{springer-LAG}.
   
  Denote by $\{\alpha_1,\dots,\alpha_r\} \subset S$ the simple roots
  which are not in $J$.  Since $J \not = \tilde S$, we have $r \ge
  1$.  Write $\varpi_i^\vee = \varpi_{\alpha_i}^\vee$.  Choose $\ell$
  a prime number different from $p$, and let $s \in \TTT$ be the image
  of
  \begin{equation*}
    \tilde s =   \dfrac{\ell - (a_2 + \cdots + a_r)}{a_1\ell} \varpi_1^\vee + 
    \dfrac{1}{\ell} \sum_{i=2}^r  \varpi_i^\vee \in V.
  \end{equation*} 
  We have written $a_i$ for the coefficient $a_{\alpha_i}$; see eq.
  \eqref{eq:high-root}.  If $r>1$, the order of $s$ is divisible by
  $\ell$ and divides $a_1 \ell$; if $J = \tilde S \setminus
  \{\alpha_1\}$, $s$ has order $a_1$.  Since $p$ is good, the order of
  $s$ is thus relatively prime to $p$.  If $\ell$ is chosen
  sufficiently large, we have $J = \{\beta \in \tilde S \mid \langle
  \beta,\tilde s \rangle \in \Z\}.$ Since $\tilde s$ lies in the
  fundamental alcove $A_o$, (the proof of) Lemma \ref{lem:R_t=R_J}
  implies that $R_s = R_J$.  Choose an element $t \in T$
  corresponding to $s \in \TTT$ as in Lemma
  \ref{lem:springer-dual-lemma}(b). By Lemma \ref{lem:ss-centralizer},
  $C_G^o(t)$ is generated by $T$ and the $\XX_\alpha$ with $\alpha \in
  R_s$; thus $C_G^o(t) = L_J$ as desired.
\end{proof}

\section{Explicit descriptions for simple and adjoint $G$}
\label{sec:explicit}

In this section, we consider $G$ simple of adjoint type.  Thus
the roots $R$ span the weight lattice $X$ over $\Z$ and
the root system is irreducible.  The
characteristic of $k$ is assumed to be good for $G$ throughout.

 The results of the previous section show that in good characteristic,
 a pseudo-Levi subgroup in the sense of this paper (connected
 centralizer of a semisimple element) is the same as a pseudo-Levi
 subgroup in the sense of \cite{sommers-generalized-bala-carter}
 (subgroup conjugate to some $L_J$).

\begin{lem}
    \label{lem:adjoint-lem-1} Let $L_J$ be a pseudo-Levi subgroup with
center $Z$.
\begin{enumerate}
\item Put $d_J = \gcd(a_\alpha \mid \alpha \in \tilde S \setminus
J)$. Then $Z/Z^o$ is cyclic of order $d_J$.
\item Every element of the character group of $Z/Z^o$ can be
represented by a root in $R$.  \end{enumerate}
\end{lem}

\begin{proof}
  Since $p$ is good, $\Z R/ \Z J$ has no $p$-torsion.  Thus the
  character group $X(Z/Z^o)$ is isomorphic to the torsion subgroup of
  $\Z R / \Z J$ as in \cite[\S 2]{sommers-generalized-bala-carter}, so
  the proof of (1) in \emph{loc. cit.} remains valid over $k$.
  
  It is also true that $X(Z/Z^o)$ is naturally isomorphic to $\Z
  \bar{R}_J / \Z J$ where $\bar{R}_J$ denotes the rational closure of
  $R_J$ in $R$.  We will show that the set $\bar{R}_J$ surjects onto
  the latter cyclic group.  Now $\bar{R}_J$ is the root system of a
  Levi subgroup of $G$, and so it contains at most one irreducible
  component of type different than type $A$.  Since the rank of $\Z
  R_J$ equals the rank of $\Z \bar{R}_J$, the type $A$ components of
  $\bar{R}_J$ play no role (every root sub-system is rationally closed
  in type $A$), and so we may assume that $\bar{R}_J$ is irreducible.
  Then $R_J$ is a  root system  with Dynkin diagram obtained by
  removing  one simple root $\alpha$ from the extended Dynkin diagram
  of $\bar{R}_J$.  Since there exists a positive root in $\bar{R}_J$
  whose coefficient on $\alpha$ is any number between $1$ and $d_J$
  (note that $d_J$ is necessarily the coefficient of the highest root
  of $\bar{R}_J$ on $\alpha$), (2) follows.
\end{proof}

\begin{lem}
  \label{lem:adjoint-lem-2}
  Let $L$ be a pseudo-Levi subgroup with center $Z$.
  \begin{enumerate}
  \item For $t \in Z$, we have 
        $L = C^o_G(tZ^o)$ if and only if $tZ^o$ generates $Z/Z^o$.
  \item If $u \in L$ is a distinguished unipotent element, then the
    group $N_G(L) \cap C_G(u)$ acts transitively on the generators of the
    cyclic group $Z/Z^o$.
  \end{enumerate}
\end{lem}

\begin{proof}
  We always have $L=C^o_G(Z) \subset C^o_G(tZ^o) \subset C^o_G(Z^o)$.
  Hence if $tZ^o$ generates $Z/Z^o$, then clearly $L=C^o_G(tZ^o)$.
  For the converse, we may assume that $L=C^o_G(t)$ by Proposition
  \ref{prop:misc-plevi}.  If $tZ^o$ fails to generate $Z/Z^o$,
  then by the previous lemma there exists a root $\beta \in R$ such
  that $\beta(t)=1$, but $\beta$ is non-trivial on $Z$.  By 
  Lemma \ref{lem:ss-centralizer} this
  contradicts the fact that $C^o_G(Z)=C^o_G(t)$, and (1) follows.

    Assertion (2) follows from \cite[Prop.
  8]{sommers-generalized-bala-carter}.
\end{proof}

\begin{prop}
  \label{prop:main-adjoint-result}
  \cite{sommers-generalized-bala-carter} 
        To a pair $(L,u)$ of a pseudo-Levi subgroup $L$ with center
  $Z$ and distinguished unipotent element $u \in L$, assign the pair
  $(u,c)$ where $c \in A(u)$ is the image of any generator of $Z/Z^o$.
  Then this map defines a bijection between the $G$-orbits on the
  pairs $(L,u)$ and the $G$-orbits on the pairs $(u,c)$.
\end{prop}

\begin{proof}
  In view of Lemma \ref{lem:adjoint-lem-2}, this follows from Theorem
  \ref{theorem:main-theorem}.
\end{proof}

To determine the isomorphism type of the groups $A(u)$ we need to
argue that the calculations in \cite{sommers-generalized-bala-carter}
remain valid over $k$.

Let $\hat G$ be the group over $\C$ with the same root datum as $G$.
Since the characteristic is good, the Bala-Carter-Pommerening theorem
shows that unipotent classes of $G$ and of $\hat G$ are parametrized
by their labeled diagram; cf. \cite[4.7 and 4.13]{jantzen:Nilpotent}.
It follows immediately that the $G$-orbits of pairs $(L,u)$ as in the
previous proposition are parametrized by the same combinatorial data
as for $\hat G$; namely, $(L,u)$ corresponds to the pair $(J,D_J)$ where
$J$ is a proper subset of $\tilde S$ and $L$ is conjugate to $L_J$
(see Proposition \ref{prop:semisimple-centralizer}), and where $D_J$ is the
labeled Dynkin diagram of the class of $u$ in $L$.  
As in the remarks
preceding \cite[Remark 6]{sommers-generalized-bala-carter}, the
$G$-orbit of $(L,u)$ identifies with the $W$-orbit of $(J,D_J)$.

Now given a unipotent class in $G$ with labeled diagram $D$, we are
left with the task of determining which pairs $(L,u)$ (up to
$G$-conjugacy) as in the previous proposition are such that $u$ has
diagram $D$ in $G$.  Since an associated cocharacter of $u$ in $L$ is
associated to $u$ in $G$ by Proposition \ref{prop:assoc-in-PL}, we may
begin with the labeled diagram of $u$ for $L$ and produce by
$W$-conjugation the labeled diagram of $u$ for $G$; see
\cite[\S3.3]{sommers-generalized-bala-carter}.  It is now clear that
our task is combinatorial: for a fixed $J \subsetneq \tilde S$, we
must find all ``distinguished'' labeled diagrams for $L_J$ which have
$D$ as a $W$-conjugate.  The calculations are carried out in
\cite[\S3.3, 3.4, 3.5]{sommers-generalized-bala-carter}, and they
remain valid for $k$.  Thanks to Proposition
\ref{prop:main-adjoint-result}, this gives a bijection between the
conjugacy classes of $A(u)$ and those of $A(\hat u)$.

According to Lemma \ref{lem:adjoint-lem-1}(1), the order of a
representative element in $A(u)$ for the class determined by the pair
$(L,u)$ is independent of the ground field.  According to
\cite[\S3.4,\S3.5]{sommers-generalized-bala-carter}, knowledge of the
conjugacy classes and the orders of representing elements in $A(\hat
u)$ are sufficient to determine the group structure. The same then
holds for $A(u)$, and we have proved:

\begin{theorem}
  \label{theorem:adjoint}
        For each unipotent element $u \in G$, let $\hat u \in \hat G$
        be a unipotent element with the same labeled diagram as $u$.
        Then $A(u) \iso A(\hat u)$. 
\end{theorem}

\providecommand{\bysame}{\leavevmode\hbox to3em{\hrulefill}\thinspace}
\providecommand{\MR}{\relax\ifhmode\unskip\space\fi MR }
\providecommand{\MRhref}[2]{%
  \href{http://www.ams.org/mathscinet-getitem?mr=#1}{#2}
}
\providecommand{\href}[2]{#2}

\end{document}